\documentclass{amsart}

\usepackage{leqno}
\usepackage{graphics}
\usepackage{amssymb}
\usepackage{amsmath}
\usepackage{amsthm}
\usepackage[latin1]{inputenc}
\usepackage{euscript}

\usepackage{color}

\newtheorem{theorem}{Theorem}[section]
\newtheorem{lemma}[theorem]{Lemma}

\theoremstyle{definition}

\newtheorem{proposition}[theorem]{Proposition}
\theoremstyle{remark}
\newtheorem{remark}[theorem]{Remark}

\numberwithin{equation}{section}

 \numberwithin{equation}{section}

\newcommand{\beq}{\begin{equation}}
\newcommand{\eeq}[1]{\label{#1}\end{equation}}

\def\O{{\Omega}}

\def\eps{{\epsilon}}

\def\j{{\mathcal{J}}}
\def\l{{\mathcal{L}}}
\def\m{{\mathcal{M}}}

\def\L{{\mathcal{L}}}

\def\M{{\mathcal{M}}}

\def\R{{\mathbb{R}}}

\def\Z{{\mathbb{Z}}}
\def\N{{\mathbb{N}}}

\newcommand{\opm}[1]{\M[{#1}]}

\newcommand{\fdem}{\vskip 0.2 pt \hfill  $\square$  }

\newtheorem{theo}{\textbf{Theorem}}[section]

\newtheorem{cla}[theo]{\textbf{Claim}}

\def\tilde{\widetilde}

\usepackage{exscale}
\renewcommand{\ldots}{\ensuremath{\dotsc}}

\newcommand{\betless}{\noalign{\vskip3pt plus 3pt minus 1pt}}
\newcommand{\bet}{\noalign{\vskip6pt plus 3pt minus 1pt}}

\def\ed{\end{document}}
\title{Existence and uniqueness
of solutions to a non-local equation with monostable nonlinearity}

\author{J\'er\^ome Coville$^{2,3}$} \author{ Juan D\'avila$^{1,2}$ \and Salom\'e Mart\'\i nez$^{1,2}$}
\address{$^1$  Departamento de Ingenier\'ia Matem\'atica\\
    Universidad de Chile
	Blanco Encalada - 2120, 5 PISO\\ 
   Santiago - Chile}
\address{$^2$   Centro de Modelamiento Matem\'atico\\
UMI 2807 CNRS-Universidad de Chile\\
Av. Blanco Encalada 2120, Piso 7\\
Santiago - Chile}
\email{jdavila@dim.uchile.cl}
\email{samartin@dim.uchile.cl}
\curraddr{$^3$  Max Planck Institut for Mathematics in the Sciences\\Inselstrasse 22\\
D-04103     Leipzig\\
Deutschland}
\email{coville@mis.mpg.de}

\date{October 31, 2006}

\begin{document}

\begin{abstract}
Let $J \in C(\R)$, $J\ge 0$, $\int_\R J = 1$ and consider the
nonlocal diffusion operator $\M[u] = J \star u - u$. We study the
equation
$$
\M u + f(x,u) = 0, \quad u \ge 0 \quad \hbox{in $\R$},
$$
where $f$ is a KPP type non-linearity, periodic in $x$. We show
that the principal eigenvalue of the linearization around zero is
well defined and a that a nontrivial solution of the nonlinear
problem exists if and only if this eigenvalue is negative. We
prove that if, additionally, $J$ is symmetric then the non-trivial
solution is unique.
\end{abstract}

\maketitle

\section{Introduction}
Reaction-diffusion equations have been used to describe a variety
of phenomena in combustion theory, bacterial growth, nerve
propagation, epidemiology, and spatial
ecology~\cite{Fisher,Fife,Kanel,Murray}. However, in many
situations, such as in population ecology, dispersal is better
described as a long range process rather than as a local one, and
integral operators appear as a natural choice. Let us mention in
particular the seminal work of Kolmogorov, Petrovsky, and Piskunov
\cite{KPP}, who in 1937 introduced a model for the
dispersion of gene fractions involving a nonlocal linear operator
and a nonlinearity of the form $u(1-u)$, which many authors now 
call a KPP-type nonlinearity.

Nonlocal dispersal operators usually take the form $\M[ u] =
\int_{\R^N}k(x,y)u(y) dy-u(x)$, where $k \ge 0$ and $\int_{\R^N}
k(y,x)dy=1$ for all $x\in \R^N$. They have been mainly used in
discrete time models \cite{KLD}, while continuous time versions
have also been recently considered in population dynamics
\cite{HMMV,LPL}. Steady state and travelling wave solutions for
single equations have been studied in the case $k(x,y)=J(x-y)$,
with $J$ even, for some specific reaction nonlinearities in
\cite{BFRW,CD,CR,BHR,CC,Schumacher}.

In this work we restrict ourselves to one dimension and take
$$
k(x,y) = J(x-y).
$$

We are interested in the existence/nonexistence and uniqueness of
solutions of the following problem:
\begin{align}
&\opm{u}+f(x,u)=0 \quad \text{ in } \R,\label{cdm.eq.noloc}
\end{align}
where $f(x,u)$ is {a KPP-type nonlinearity, periodic in $x$,} and
\begin{align}
\label{operator M} \opm{u}:=J\star u -u.
\end{align}
We assume that $J$ satisfies
\begin{align}
\label{J1} & J \in C(\R), \quad J\ge 0, \quad {\int_\R J = 1} ,
\\
\bet
\label{J2} & \hbox{there exist $a<0<b$ such that $J(a)>0$,
$J(b)>0$.}
\end{align}

On $f$ we assume that
\begin{align}
\label{hyp f1} \left\{
\begin{aligned}
& \hbox{$f \in C(\R\times[0,\infty))$ and is differentiable with respect to $u$,} \\
\betless
& \hbox{for each $u$, $f(\cdot,u)$ is periodic with period $2 R$,}
\\
\betless
& \hbox{$f_u(\cdot,0)$ is Lipschitz,} \\
\betless
& \hbox{$f(\cdot,0)\equiv0$ and $f(x,u)/u$ is decreasing with
respect to $u$,} \\
\betless
& \hbox{there exists $M>0$ such that $f(x,u)\le0$ for all $u \ge
M$ and all $x$.}
\end{aligned}
\right.
\end{align}
The model example of such a nonlinearity is
$$
f(x,u) = u ( a(x) - u ),
$$
where $a(x)$ is periodic and Lipschitz.

In a recent work, Berestycki, Hamel, and Roques \cite{BHR} studied
the analogue of \eqref{cdm.eq.noloc} with a divergence operator in a
periodic setting. More precisely, they considered
\begin{align}
\label{local eq} - \nabla \cdot (A(x) \nabla u ) = f(x,u), \quad
x\in \R^N, \quad u \ge 0,
\end{align}
where $A(x)$ is a symmetric matrix of class $C^{1,\alpha}$,
periodic with respect to all variables and uniformly elliptic, and
$f$ is $C^1$ and satisfies \eqref{hyp f1}. They showed existence
of nontrivial solutions provided the linearization of the equation
around zero has a negative first periodic eigenvalue.

We prove the following result.

\begin{theorem}
\label{thm0}
Assume $J$ satisfies \eqref{J1}, \eqref{J2} and $f$
satisfies \eqref{hyp f1}. Then there exists a nontrivial, periodic solution of \eqref{cdm.eq.noloc} if and only
if
$$
\lambda_1(\m+f_u(x,0))<0,
$$
where $\lambda_1$ is the principal eigenvalue of the linear
operator {$-(\m+f_u(x,0))$} in the set of $2R$-periodic continuous
functions. Moreover, if $\lambda_1 \ge 0$, then any nonnegative
bounded solution is identically zero.
\end{theorem}

To prove Theorem~\ref{thm0}, we first need to show that the
principal periodic eigenvalue of $-(\m+f_u(x,0))$ is well defined.
Let us introduce some notation:
\begin{align*}
& C_{per}(\R) = \{ u:\R\to\R \; | \; \hbox{$u$ is continuous and $2R$-periodic} \},\\
\bet
& C^{0,1}_{per}(\R) = \{ u:\R\to\R \; | \; \hbox{$u$ is Lipschitz and $2R$-periodic} \}.
\end{align*}

\begin{theorem}
\label{cdm.th.pev}
Suppose $a(x)\in C^{0,1}_{per}(\R)$. Then the operator $- (\M +
a(x))$ has a unique principal eigenvalue $\lambda_1$ in
$C_{per}(\R)$; that is, there is a unique $\lambda_1 \in \R$ such
that
\begin{align}
\label{eigenfunction} \opm{\phi_1}+a(x)\phi_{1}=-\lambda_1\phi_1
\quad \hbox{in $\R$}
\end{align}
admits a positive solution $\phi_1 \in C_{per}(\R)$. Moreover,
$\lambda_1$ is simple, that is, the space of $C_{per}(\R)$
solutions to \eqref{eigenfunction} is one dimensional.
\end{theorem}

In \cite{BHR} the authors proved that \eqref{local eq} has at most
one nontrivial bounded solution, and that it has to be periodic. A
similar result is true for the nonlocal problem
\eqref{cdm.eq.noloc}, but this time we need $J$ to be symmetric,
that is,
\begin{align}
\label{J symmetric} J(x)=J(-x) \quad \hbox{for all }x\in\R.
\end{align}
Note, however, that for the existence result, Theorem~\ref{thm0},
we do not need this \hbox{condition.}

\begin{theorem}
\label{uniqueness J symmetric}
Assume $J$ satisfies \eqref{J1},
\eqref{J2}, \eqref{J symmetric} and $f$ satisfies \eqref{hyp f1}. Let
$u$ be a nonnegative, bounded solution
to \eqref{cdm.eq.noloc} and let $\lambda_1$ be the principal
eigenvalue of the operator $-(\M + f_u(x,0))$ with periodic
boundary conditions.

{\rm (a)} If $\lambda_1<0$, then either $u\equiv 0$ or $u\equiv p$, where
$p$ is the positive periodic solution of Theorem~$\ref{thm0}$.

{\rm (b)} If $\lambda_1 \ge 0$, then $u\equiv 0$.
\end{theorem}

Part (b) of the preceding theorem is already covered in
Theorem~\ref{thm0} and does not depend on the symmetry of $J$.

When $f $ is independent of $x$ and satisfies \eqref{hyp f1}, 
the principal eigenvalue of $-( \M + f'(0) )$ is given by
$\lambda_1 = - f'(0)$ and $\phi_1$ is just a constant. Thus in
this case Theorem~\ref{thm0} says that a bounded, nonnegative,
nontrivial solution exists if and only if $f'(0)>0$, and this
solution is just the constant $u_0$ such that $f(u_0)=0$. Assuming
that $J$ is symmetric, Theorem~\ref{uniqueness J symmetric} then
implies that the constant $u_0$ is the unique solution in the
class of nonnegative, bounded functions.

Recently, considering a nonperiodic nonlinearity $f$, Berestycki,
Hamel, and Rossi \cite{BHRo} analyzed the analogue of
Theorem~\ref{uniqueness J symmetric} for general elliptic
operators in $\R^N$, finding sufficient conditions that ensure
existence and uniqueness of a positive bounded solution. It is
natural to ask whether the periodicity of $f$ and the symmetry of
$J$ are crucial hypotheses in Theorem~\ref{uniqueness J
symmetric}. We believe that this is the case, since a general
nonlocal operator such as \eqref{operator M} may contain a transport
term, and a standing wave connecting the steady states of the
system could appear. We shall investigate further this issue in a
forthcoming work.

Hypothesis \eqref{J2} implies that the operator $\M$ satisfies the
strong maximum principle. Suppose, for instance, that $J$
satisfies \eqref{J1}, \eqref{J2}. If $u \in C(\R)$ satisfies
$\M[u] \ge 0 $ in $\R$, then $u$ cannot achieve a global maximum
without being constant (see \cite{C}). However, we will need the
following version.

\begin{theorem}
\label{strong max p a.e.}
Assume $J$ satisfies \eqref{J1},
\eqref{J2} and let $c \in L^\infty(\R)$. If $u \in L^\infty(\R)$
satisfies $u \le 0 $ a.e.\ and $\M[u] + c(x) u \ge 0$ a.e.\ in
$\R$, then ${\rm ess\ sup}_K u <0 $ for all compact $K\subset \R$ or $u
=0 $ a.e.\@ in $\R$.
\end{theorem}

If $f$ satisfies the stronger hypothesis that, for any $x$, $f(x,u)$ is
concave with respect to $u$, then actually the periodic solution
$p$ of Theorem~\ref{thm0} is continuous. To see this notice that
from the strong maximum principle, Theorem~\ref{strong max p
a.e.}, $J \star p >0$ in $\R$. The concavity of $f$ with respect
to $u$ implies that for any $x$ the map $ u \mapsto u - f(x,u)$ is
strictly increasing whenever $u - f(x,u)>0$. Then from the
continuity of $J \star p$ and \eqref{cdm.eq.noloc},
which can be rewritten as in the form $J \star p = p - f(x,p)$, we
deduce that $p$ is continuous.

In section~\ref{spectral theory} we review some spectral theory
and give the argument of Theorem~\ref{cdm.th.pev}. Then we prove
Theorem~\ref{thm0} in section~\ref{section Existence of solutions}
and the uniqueness result, Theorem~\ref{uniqueness J symmetric}(a),
in section~\ref{Uniqueness when J is symmetric}. We leave for an
appendix a proof of Theorem~\ref{strong max p a.e.}.

\section{Some spectral theory}
\label{spectral theory}

In this section we deal with the principal eigenvalue problem
\eqref{eigenfunction}. Before stating our result, let us recall
some basic spectral results for positive
operators due to Edmunds, Potter, and Stuart \cite{EPS} which are
extensions of the Krein--Rutmann theorem for positive noncompact
operators.

A cone in a real Banach space $X$ is a nonempty closed set $K$
such that for all $x, y \in K$ and all $\alpha \ge 0$ one has $x +
\alpha y \in K$, and if $x\in K$, $-x \in K$, then $x=0$. A cone
$K$ is called reproducing if $X = K-K$. A cone $K$ induces a
partial ordering in $X$ by the relation $x\le y$ if and only if
$x-y \in K$. A linear map or operator $T:X\to X$ is called
positive if $T(K)\subseteq K$. The dual cone $K^*$ is the set of
functionals $x^* \in X^*$ which are positive, that is, such that
$x^*(K) \subset [0,\infty)$.

If $T:X\to X$ is a bounded linear map on a complex Banach space X,
its essential spectrum (according to Browder \cite{browder})
consists of those $\lambda$ in the spectrum of $T$ such that at
least one of the following conditions holds: (1)~the range of
$\lambda I - T$ is not closed, (2)~$\lambda$ is a limit point of
the spectrum of $T$, (3)~$\cup_{n=1}^\infty \ker( (\lambda I -
T)^n)$ is infinite dimensional. The radius of the essential
spectrum of $T$, denoted by $r_e(T)$, is the largest value of
$|\lambda|$ with $\lambda$ in the essential spectrum of $T$. For
more properties of $r_e(T)$ see \cite{nussbaum}.

\begin{theorem}[Edmunds, Potter, and Stuart~\cite{EPS}]
\label{cdm.th.eps}
Let K be a reproducing cone in a real Banach space X, and let
$T\in\L(X)$ be a positive operator such that $T^p(u)\ge cu$ for
some $u\in K$ with $\|u\|=1$, some positive integer $p$, and some
positive number $c$. Then if $c^{\frac{1}{p}}>r_e(T)$, $T$ has an
eigenvector $v\in K$ with associated eigenvalue $\rho\ge
c^{\frac{1}{p}}$ and $T^*$ has an eigenvector $v^*\in K^*$
corresponding to the eigenvalue $\rho$.
\end{theorem}

A proof of this theorem can be found in \cite{EPS}. If the cone
$K$ has nonempty interior and $T$ is strongly positive, i.e.,
$u\ge0$, $u\not=0$ implies $T u \in {\rm int}(K)$, then $\rho$ is the
unique $\lambda \in \R$ for which there exists nontrivial $v \in K$
such that $T v = \lambda v$ and $\rho$ is simple; see~\cite{zeidler}.

{\it Proof of Theorem}~\ref{cdm.th.pev}.
For convenience, in this
proof we write the eigenvalue problem
\begin{align*}
\M[u] + a(x)u = -\lambda u
\end{align*}
in the form
\begin{align}
\label{cdm.eq.pev2}
\L[u] + b(x) u = \mu u,
\end{align}
where
$$
\L[u] = J\star u, \quad b(x) = a(x) + k, \quad \mu = -\lambda + 1 + k,
$$
and $k>0$ is a constant such that $\inf_{[-R,R]} b >0$.

Observe that $\L:C_{per(\R)} \to C_{per}(\R)$ is compact
($C_{per}(\R)$ is endowed with the norm
$\|u\|_{L^\infty([-R,R])}$). Indeed, let $u_n \in C_{per}(\R)$ be
a bounded sequence, say\linebreak
$\|u_n \|_{L^\infty([-R,R])} \le B$. Let
$\eps>0$ and let $A$ be large enough so that $\int_{|x|\ge A} J \le
\eps$. Since $J$ is uniformly continuous in $[-R-2A,R+2A]$ there
is $\delta>0$ such that $| J(z_1 ) - J(z_2)| \le
\frac{\eps}{2(A+R)}$ for $z_1,z_2 \in [-R-2A,R+2A]$ with $|z_1
-z_2|\le \delta$. Then for $x_1,x_2 \in [-R,R]$,
\begin{align*}
|\L[u_n](x_1) - \L[u_n](x_2) | & \le \int_\R | J(x_1 - y ) -
J(x_2-y) | \, | u_n(y)| \, d y
\\
\bet
& \le 2 B \eps + B \int_{-R-A}^{R+A} | J(x_1 - y ) - J(x_2-y) |
\, d y
\\
\bet
& \le 3 B \eps.
\end{align*}
This shows that $\L[u_n]$ is equicontinuous, and therefore by the
Arzel\`a--Ascoli theorem, $ \L[u_n]$ is relatively compact.

Let us now establish some useful lemma.

\begin{lemma}\label{cdm.lem.pev}
Suppose $b(x)\in C^{0,1}(\R)$ is $2R$-periodic, $b(x)>0$, and let
$\sigma:=\max_{ [-R,R]} b(x)$. Then there exist $p \in \N, \delta>0$,
and $u \in C_{per}(\R)$, $u \ge 0$, $u \not \equiv 0$, such that
$$
\l^{p}u+b(x)^p u\ge (\sigma^p+\delta) u.
$$
\end{lemma}
\unskip

Observe that the proof of Theorem~\ref{cdm.th.pev} will then
easily follow from the above lemma. Indeed, if the lemma holds,
then since $u$ and $b$ are nonnegative and $\L$ is a positive
operator, we easily see that
$$
(\l+b(x))^{p}[u] \ge \l^{p}[u]+b(x)^{p}u \ge (\sigma^p+\delta)u .
$$
Using the compactness of the operator $\l$, we have
$r_e(\l+b(x))=r_e(b(x)) {=} \sigma$, and thus
$(\sigma^p+\delta)^{\frac{1}{p}}>r_e(\l+b(x))$ and Theorem~\ref{cdm.th.eps} applies.
Finally, we observe that the principal eigenvalue is simple since
the cone of positive $2R$-periodic functions has nonempty
interior and, for a sufficiently large $p$, the operator $(\l +
b)^p $ is strongly positive.
\fdem

Let us now turn our
attention to the proof of the above lemma.

{\it Proof of Lemma~$2.2$}. Recall that for $p\in { \N \setminus
\{0\} }$, $J\star^{p} u:=J\star(J\star^{p-1}u)$ is well defined by
induction and satisfies $J\star^{p}u=\j_p\star u$ with $\j_p$
defined as follows:
$$
\begin{array}{lc}
\j_p:=&\underbrace{J\star J\star \cdots \star J\star J}\\
 & p \text{ times }
\end{array}.
$$

By \eqref{J2} it follows that there exists $p\in \N$ such that
$\inf_{(-2R-1,2R+1)} \j_{p} > 0 $.
Using the definition of $\l$, a short computation shows that
$$
\l^p[u]:= \int_{-R}^{R}\tilde \j_p(x,y)u(y)\,dy
$$
with $\tilde\j_p(x,y)=\sum_{k\in\Z}\j_p(x+{2kR}-y)$. Following the idea
of Hutson et~al.\ \cite{HMMV}, consider now the following function:
{\begin{equation*} v(x):=\left\{
\begin{array}{ll}
\frac{\eta(x) }{b^p(x_0)-b^p(x) +\gamma} &\text{ in } \O_{2\eps} :=(x_0-2\eps ,x_0+2\eps),\\
\betless
0 & \text{elsewhere},
\end{array}
\right.
\end{equation*}
where $x_0\in (-R,R)$ is a point of maximum of $b(x)$, $\eps>0$ is
chosen such that $(x_0-2\eps ,x_0+2\eps)\subset(-R,R)$, $\gamma$
is a positive constant that we will define later on, and $\eta$ is
a smooth function such that $0\le \eta \le 1$, $\eta(x) = 1$ for
$|x-x_0| \le \eps$, $\eta(x) = 0$ for $|x-x_0| \ge 2\eps$.} Let us
compute $\l^p[v]+b^p(x)v-\sigma^p v$:
\begin{align*}
\l^p[v]+b^p(x)v-\sigma^p v&=\int_{x_0-\eps}^{x_0+\eps}\tilde
\j_p(x,y)\frac{dy}{b^p(x_0)-b^p(y) +\gamma} +
\int_{\O_{2\eps}\setminus\O_\eps}\tilde \j_p(x,y) v(y) \,dy
\\
\bet
& \quad +(b^p(x)-b^p(x_0))v
\\
\bet
&\ge \int_{x_0-\eps}^{x_0+\eps}
\tilde\j_p(x,y)\frac{dy}{b^p(x_0)-b^p(y) +\gamma}
+(b^p(x)-b^p(x_0))v
\\
\bet
&\ge \int_{x_0-\eps}^{x_0+\eps}
\tilde\j_p(x,y)\frac{dy}{b^p(x_0)-b^p(y) +\gamma} -1 .
\end{align*}
Using that $\inf_{(-2R-1,2R+1)} \j_{p} > 0 $, it follows that
$\tilde \j_p(x,y) \ge c >0$ for $x,y \in (-R,R)$. Hence
$$
\int_{x_0-\eps}^{x_0+\eps}\tilde
\j_p(x,y)\frac{dy}{b^p(x_0)-b^p(y) +\gamma}\ge
c\int_{x_0-\eps}^{x_0+\eps}\frac{dy}{k|x_0-y| +\gamma},
$$
{where $k$ is the Lipschitz constant for $b^p$.} Using this
inequality in the above estimate yields
$$
 \l^p[v]+b^p(x)v-\sigma^p v\ge c\int_{x_0-\eps}^{x_0+\eps}\frac{dy}{k|x_0-y| +\gamma}-1.
$$
Therefore we have
\begin{align*}
\l^p[v]+b^p(x)v-(\sigma^p+\delta) v&\ge\frac{2c}{k}\log{\left(1+\frac{k\eps}{\gamma}\right)}-1 -\delta v\\
\bet
&\ge \frac{2c}{k}\log{\left(1+\frac{k\eps}{\gamma}\right)}-1
-\frac{\delta}{\gamma} .
\end{align*}
Choosing now $\gamma>0$ small so that
$\frac{2c}{k}\log{(1+\frac{k\eps}{\gamma})}-1>\frac{1}{2}$ and
$\delta=\frac{\gamma}{4}$, we end up with
$$
\l^p[v]+b^p(x)v-(\sigma^p+\delta) v\ge \frac{1}{4}>0.
$$
\fdem

\section{Existence of solutions}
\label{section Existence of solutions}\quad

{\it Proof of Theorem}~\ref{thm0}.
We follow the {argument}
developed by Berestycki, Hamel, and Roques in \cite{BHR}.

First assume that $\lambda_1< 0$.
From Theorem~\ref{cdm.th.pev} there exists a positive
eigenfunction $\phi_1$ such that
$$
\opm{\phi_1}+f_u(x,0)\phi_1=-\lambda_1\phi_1\ge 0.
$$

Computing $\opm{\eps\phi_1} +f(x,\eps\phi_1)$, it follows that
\begin{align*}
\opm{\eps\phi_1} +f(x,\eps\phi_1)&=f(x,\eps\phi_1)-f_u(x,0)\eps\phi_1-\lambda_1 {\eps} \phi_1\\
\bet
&=-\lambda_1\eps\phi_1+o(\eps\phi_1)>0.
\end{align*}
Therefore, for $\eps>0$ small, $\eps\phi_1$ is a periodic
subsolution of \eqref{cdm.eq.noloc}. By definition of $f$, any
constant $M$ sufficiently large is a periodic supersolution of the
problem. Choosing $M$ so {large} that $\eps\phi_1\le M$ and using
a basic iterative scheme yields the existence of a positive
periodic solution $u$ of \eqref{cdm.eq.noloc}.

Let us now turn our attention to the nonexistence setting and assume that $\lambda_1\ge 0$.

\noindent
Let $u$ be a bounded nonnegative solution of
\eqref{cdm.eq.noloc}. Observe that $\gamma\phi_1$ is a periodic
supersolution for any positive $\gamma$. Indeed,
\begin{align*}
\opm{\gamma\phi_1}+f(x,\gamma\phi_1)& < \opm{\gamma\phi_1}+f_{u}(x,0)\gamma\phi_1\\
\bet
& \le -\lambda_1\gamma\phi_1\le 0.
\end{align*}
Since $\phi_1\ge\delta$ for some positive $\delta$ we may define
the following quantity:
$$
\gamma^*:=\inf\{\gamma>0|u\le\gamma\phi_1\}.
$$
We have the following claim.

\begin{cla}
$\gamma^*=0$.
\end{cla}

Observe that we end the proof of the theorem by proving the above
claim.

{\it Proof of the claim}.
Assume that $\gamma^*>0$. Since $v := u
- \gamma^* \phi_1$ satisfies $v \le 0 $ in $\R$ and
$$
\M[v] + c(x) v \ge 0 \quad \hbox{in $\R$},
$$
where $c(x) = \frac{f(x,u)- f(x, \gamma^*\phi_1)}{v}$ by the
strong maximum principle, Theorem~\ref{strong max p a.e.}, we have
the following possibilities:

\begin{itemize}
\item
either $ u \equiv \gamma^*\phi_1$, or

\item
there exists a sequence of points $(x_n)_{n\in\N}$ such that
{$|x_n|\to +\infty$} and\linebreak
$\lim_{n\to+\infty}\gamma^*\phi_1(x_n)-u(x_n)=0$.
\end{itemize}

\noindent
In the first case we get the following contradiction:
$$
0=\opm{ \gamma^* \phi_1}+f(x, \gamma^* \phi_1) < \opm{\gamma^*
\phi_1}+f_{u}(x,0)\gamma^*\phi_1\le 0.
$$
Hence $\gamma^*=0$.

{In the second case} we argue as follows. Let $(y_n)_{n\in\N}$ be
a sequence of points satisfying, for all $n$, $y_n\in [-R,R]$ and
$x_n-y_n \in 2 R\Z$. Up to extraction {of a subsequence}, $y_n \to
\bar y$. Now consider the following sequence of functions
$u_n:=u(.+x_n)$, $\phi_n:=\phi_1(.+x_n)$, and $w_n:=\gamma^*\phi_n
-u_n$ so that $w_n>0$ in $\R$. Since $\m$ is translation invariant
and $f$ is periodic, $u_n$ and $\phi_n > 0$ satisfy
\begin{align*}
&\opm{u_n}+f(x+y_n,u_n)=0 \quad \hbox{in $\R$},\\
\bet
&\opm{\gamma^*\phi_n}+f_u(x+y_n,0) \gamma^*\phi_n \le 0 \quad
\hbox{in $\R$} .
\end{align*}
It follows that
$$
J\star w_n \le a_n(x) w_n,
$$
where
$$
a_n(x) = 1-\frac{ \gamma^* f_u(x+y_n,0) \phi_n - f(x+y_n,u_n)
}{\gamma^* \phi_n - u_n} .
$$
Since $w_n>0$ we see that $a_n$ is well defined and $a_n \ge 0$. Using
that $f(x,u)/u$ is nonincreasing with respect to $u$ we have $
f(x,\gamma^* \phi_n) \le \gamma^* f_u(x,0) \phi_n$. This implies
\begin{align*}
\frac{ \gamma^* f_u(x+y_n,0) \phi_n - f(x+y_n,u_n) }{\gamma^*
\phi_n - u_n} \ge \frac{ f(x+y_n, \gamma^* \phi_n) -
f(x+y_n,u_n) }{\gamma^* \phi_n - u_n} \ge -C.
\end{align*}
Thus
$$
0 \le a_n \le C +1 \quad \hbox{ in $\R$ \quad for all $n$,}
$$
with $C$ independent of $n$. Observe that
$$
J \star w_n(0) = a_n(0)(\gamma^* \phi_1(x_n) - u(x_n) )\to 0,
$$
which implies
$$
\int_\R J(-y) w_n(y) \,dy \to 0 \quad \hbox{ as } n\to+\infty.
$$
Similarly,
$$
J\star J\star w_n (0) = J\star(a_n w_n)(0) = \int_\R J(-y) a_n(y)
w_n(y)\,dy,
$$
but
$$
\int_\R J(-y) a_n(y) w_n(y) \, dy \le \|a_n\|_{L^\infty} \int_\R
J(-y) w_n(y) \, dy \to 0.
$$
Hence
$$
J \star J \star w_n(0) = \int_\R (J\star J)(-y) w_n(y) \, d y \to
0 \quad \hbox{as $n\to+\infty$}.
$$
Defining
$$
\j_k:=\underbrace{J\star \cdots \star J}_{k \text{ times }}\ ,
$$
we see that for any fixed $k \in \N$,
$$
\int_\R \j_k (-y) w_n(y) \, d y \to 0 \quad \hbox{as
$n\to+\infty$}.
$$
By \eqref{J2} the support of $\j_k$ increases to all of $\R$ as $k
\to +\infty$. Thus we may find a new subsequence such that $w_n
\to 0$ a.e.\@ in $\R$ as $n \to+\infty$. Since $\phi_1$ is
periodic and continuous, $\phi_n(x) \to \bar \phi(x)$ uniformly
with respect to $x$, where $\bar \phi(x)= \phi(x+\bar y)$. Hence
$\bar u(x)=\lim_{n\to+\infty} u_n(x)$ exists a.e.\@ and is given
by $ \bar u(x) = \gamma^* \bar \phi$ . By dominated convergence,
$\bar u$ is a solution to
\begin{align*}
\M[\bar u] + f(x+\bar y,\bar u)=0,
\end{align*}
while by uniform convergence
\begin{align*}
\opm{\gamma^*\bar\phi}+f_u(x+\bar y,0) \gamma^*\bar \phi \le 0
\quad \hbox{in $\R$}.
\end{align*}
Since $\bar u = \gamma^*\bar \phi$ it follows that $f(x+\bar y,
\gamma^* \bar \phi) \equiv f_u(x+\bar y,0) \gamma^*\bar \phi$.
This contradicts the fact that $f(x,u)/u$ is decreasing in $u$.
Hence, $\gamma^*=0$.
\fdem 

\section{\boldmath Uniqueness when $J$ is symmetric}
\label{Uniqueness when J is symmetric}

Throughout this section we assume that $J$ is symmetric. For the
proof of Theorem~\ref{uniqueness J symmetric} we follow the ideas
in \cite{BHR}.

{\it Proof of Theorem}~\ref{uniqueness J symmetric}.
Part (b) of this theorem is contained in Theorem~\ref{thm0} so we
concentrate on part (a).

Let $p$ denote the positive periodic solution to
\eqref{cdm.eq.noloc} constructed in Theorem~\ref{thm0} and let
$u\ge0$, $u\not\equiv 0$ be a bounded solution. We will prove that
$u \equiv p$.

We show first that $u \le p$. Set
$$
\gamma^*:=\inf\{\gamma>0 \; | \; u\le\gamma p\}.
$$
Note that $\gamma^*$ is well defined because $u$ is bounded and
$p$ is bounded below by a positive constant. We claim that
$$
\gamma^* \le 1.
$$
Suppose that $\gamma^*>1$ and note that $u \le \gamma^* p$. 
By Theorem~\ref{strong max p a.e.} either $u \equiv \gamma^* p$ or ${\rm ess}\ \inf_K (\gamma^* p - u ) >0$ 
for all compact $K \subset \R$. The first possibility leads to $f(x,\gamma^* p) = \gamma^* f(x,p)$ for all $x\in\R$,
which is not possible if $\gamma^*>1$. In the second case there exists a sequence $(x_n)_{n\in\N}$ such that {$|x_n|\to
+\infty$} and $\lim_{n\to+\infty}\gamma^* p(x_n)-u(x_n)=0$. Let
$(y_n)_{n\in\N}$ be a sequence satisfying $y_n\in [-R,R]$ and
$x_n-y_n =k_n 2R$ for some $k_n \in \Z$. We may assume that $y_n
\to \bar y$. Let $u_n:=u(.+x_n)$, which satisfies
\begin{align*}
\opm{u_n}+f(x+y_n,u_n)=0.
\end{align*}
Let $w_n= \gamma^* p(.+y_n) - u_n \ge 0$. Then $w_n>0$ in $\R$ and
$$
J\star w_n = a_n(x) w_n,
$$
where
$$
a_n(x) = 1-\frac{ \gamma^* f(x+y_n,p(x+y_n)) - f(x+y_n,u_n(x))
}{\gamma^*p(x+y_n) - u_n(x)} .
$$
Since $w_n>0$ we deduce that $a_n$ is well defined and $a_n \ge 0$.
Using that $f(x,u)/u$ is nonincreasing with respect to $u$ and
the fact that $\gamma^*>1$, we have $ f(x,\gamma^* p) \le \gamma^*
f(x,p)$. This implies
\begin{align*}
\frac{ \gamma^* f(x,p) - f(x,u) }{\gamma^*p - u} \ge \frac{
f(x,\gamma^* p) - f(x,u) }{\gamma^*p - u} \ge -C.
\end{align*}
Thus
$$
0 \le a_n \le C +1 \quad \hbox{ in $\R$ \quad for all $n$,}
$$
with $C$ independent of $n$. Observe that
$$
J \star w_n(0) = a_n(0)(\gamma^* p(y_n) - u(x_n)) =
a_n(0)(\gamma^* p(x_n) - u(x_n) )\to 0,
$$
which implies
$$
\int_\R J(-y) w_n(y) \,dy \to 0 \quad \hbox{ as } n\to+\infty.
$$
Similarly,
$$
J\star J\star w_n (0) = J\star(a_n w_n)(0) = \int_\R J(-y) a_n(y)
w_n(y)\,dy,
$$
but
$$
\int_\R J(-y) a_n(y) w_n(y) \, dy \le \|a_n\|_{L^\infty} \int_\R
J(-y) w_n(y) \, dy \to 0.
$$
Hence
$$
J \star J \star w_n(0) = \int_\R (J\star J)(-y) w_n(y) \, d y \to
0 \quad \hbox{as $n\to+\infty$}.
$$
Defining
$$
\j_k:=\underbrace{J\star \cdots \star J}_{k \text{ times }}\ ,
$$
we see that for all $k \in \N$,
$$
\int_\R \j_k (-y) w_n(y) \, d y \to 0 \quad \hbox{as
$n\to+\infty$}.
$$
Hypothesis \eqref{J2} implies that the support of $\j_k$ converges
to all of $\R$ as $k \to +\infty$. Therefore, for a subsequence,
$w_n \to 0$ a.e.\@ in $\R$ as $n \to+\infty$. Since $p$ is
periodic, for possibly a new subsequence $p(x + y_n) \to p(x +\bar
y)$ a.e. Hence, $\bar u(x)=\lim_{n\to+\infty} u_n(x)$ exists
a.e.\@ and by dominated convergence, $\bar u$ is a solution to
\begin{align}
\label{limit shifted} \M[\bar u] + f(x+\bar y,\bar u)=0.
\end{align}
But since $w_n \to 0$ a.e.\@ we have $\bar u = \gamma^* p(\cdot +
\bar y)$. Thus $\gamma^* p(\cdot+\bar y)$ is a solution to
\eqref{limit shifted}, which is impossible for $\gamma^*>1$ as
argued before.

The proof that $p \le u$ is analogous, but a key point is to prove
first that under the conditions of Theorem~\ref{uniqueness J
symmetric} any nontrivial, nonnegative solution is bounded below
by a positive constant. This is the content of
Proposition~\ref{positive lower bound}.
\fdem 

\begin{proposition}
\label{positive lower bound}
Assume that $J$ satisfies \eqref{J1},
\eqref{J2}, and \eqref{J symmetric}, $f$ satisfies \eqref{hyp f1},
and that the operator $-(\M - f_u(x,0))$ has a negative principal
periodic eigenvalue. Suppose that $u$ is a nonnegative, bounded
solution to \eqref{cdm.eq.noloc}. Then $u\equiv 0$ or there exists
a constant $c>0$ such that
$$
u(x) \ge c \quad \hbox{for all } x\in \R.
$$
\end{proposition}
\unskip

The basic tool to prove Proposition~\ref{positive lower bound},
following an idea in \cite{BHR}, is to study the principal
eigenvalue of the linearized operator in bounded domains. More
precisely, let $ \Omega = (-r,+r)$ and $a \colon\Omega \to \R$ be
Lipschitz. We consider the eigenvalue problem in $\Omega$ with
``Dirichlet boundary condition'' in the following sense:
\begin{align}
\label{truncated eigenvalue} \left\{
\begin{aligned}
&\M[\varphi] + a(x) \varphi = - \lambda \varphi \quad \hbox{in }
\Omega,
\\
\betless
& \varphi(x) = 0 \quad \hbox{for all } x \not\in \Omega ,\\
\betless
&\hbox{$\varphi|_{\overline \Omega}$ is continuous}.
\end{aligned}
\right.
\end{align}

We show that the principal eigenvalue for \eqref{truncated
eigenvalue} exists and converges to the principal periodic
eigenvalue as $r\to+\infty$. The first step is to establish
variational characterizations of these eigenvalues, which is the
argument that requires the symmetry of $J$.

\begin{lemma}
\label{principal eigenvalue omega}
Let $\Omega \subset \R$ be a
bounded open interval. Assume that $J$ satisfies \eqref{J1},
\eqref{J2}, and \eqref{J symmetric}, and let $a \colon\Omega \to \R$
be Lipschitz. Then there exists a smallest $\lambda_1$ such that
\eqref{truncated eigenvalue} has a nontrivial solution. This
eigenvalue is simple and the eigenfunctions are of constant sign
in $\Omega$. Moreover,
\begin{align}
\label{variational eigenvalue Omega}
\lambda_1 = \min_{\varphi \in
C(\overline \Omega) } - \frac{\int_{\Omega} (\M[\tilde \varphi] +
a(x) \varphi) \varphi }{\int_\Omega \varphi^2},
\end{align}
where $\tilde \varphi$ denotes the extension by $0$ of $\varphi$ to
$\R$ and the minimum is attained.
\end{lemma}

The statement and the proof are analogous to those of Theorem~3.1 in
\cite{HMMV} except that here we do not assume that $J(0)>0$. A
different formula for the principal eigenvalue with a Dirichlet
boundary condition appears in \cite{chasseigne-chaves-rossi},
where it is used to characterize the rate of decay of solutions to
a linear evolution equation.

\begin{proof}
Define the operator $X [\varphi] = \int_\Omega J(x-y) \,
\varphi(y) \, d y $ for $\varphi \in C(\overline \Omega)$. Then
$X :C(\overline \Omega) \to C(\overline \Omega)$ is compact. Let
$c_0>0$ be such that $ \inf_{\Omega} a(x) + c_0 > 0 $ and define
$\tilde a = a + c_0$. The eigenvalue problem \eqref{truncated
eigenvalue} is equivalent to the following: find $\varphi \in C(\overline
\Omega)$ and $\lambda \in \R$ such that
\begin{align*}
X[\varphi] + \tilde a \varphi = ( - \lambda + 1 + c_0)\varphi
\quad \hbox{in $\Omega$}.
\end{align*}
A calculation similar to Lemma~\ref{cdm.lem.pev} shows that there
exists an integer $p$, $u \in C( \overline \Omega)$, and $\delta >
0 $ such that
\begin{align}
\label{bound sigma0} (X + \tilde a)^p u \ge \left( \left(\max_{\overline
\Omega} \tilde a\right)^p + \delta \right) u \quad \hbox{in $\Omega$}.
\end{align}
Using Theorem~\ref{cdm.th.eps} we deduce that the operator $X +
\tilde a$ has a unique principal eigenvalue $\rho>0$ and a
principal eigenvector $\varphi_1 \in C(\overline \Omega)$. Let
$\lambda = 1+c_0-\rho$ so that $X[\varphi_1] + a(x) \varphi_1 =
(1- \lambda) \varphi_1$. From \eqref{bound sigma0} we deduce that
$\sigma_+ $ defined by
\begin{align}
\label{sigma +} \sigma_+ = \sup_{\varphi \in C(\overline \Omega) }
\frac{\int_{\Omega} (X[ \varphi] + a(x) \varphi) \varphi
}{\int_\Omega \varphi^2}
\end{align}
satisfies
\begin{align}
\label{bound sigma+} \sigma_+ \ge 1-\lambda > \max_{\overline
\Omega} a.
\end{align}
Now, using the same argument as in \cite{HMMV} we deduce that
the supremum in \eqref{sigma +} is achieved. Indeed, it is standard
\cite{brezis} that the spectrum of $\hat{X} + a(x)$ is to the left of
$\sigma_+$ and that there exists a sequence $\varphi_n \in
C(\overline \Omega)$ such that $\|\varphi_n \|_{L^2(\Omega)}=1$
and $ \| (X + a(x) - \sigma_+) \varphi_n \|_{L^2(\Omega)} \to 0$
as $n\to+\infty$. By compactness of $X:L^2(\Omega) \to
C(\overline\Omega)$ for a subsequence, $\lim_{n\to+\infty}
X[\varphi_n] $ exists in $C(\overline \Omega)$. Then, using
\eqref{bound sigma+}, we see that $\varphi_n \to \varphi$ in
$L^2(\Omega)$ for some $\varphi $ and $(X+a)\varphi = \sigma_+
\varphi$. This equation implies $ \varphi \in C(\overline \Omega)$,
and hence $\sigma_+ $ is a principal eigenvalue for the operator
$X$ and by uniqueness of this eigenvalue we have $\sigma_+ = 1 -
\lambda$.
\qquad\end{proof}

\begin{lemma}
\label{lema formulas eig R}
Assume that $J$ satisfies \eqref{J1},
\eqref{J2}, and \eqref{J symmetric} and that $a\colon\R\to\R$ is a
$2R$-periodic, Lipschitz function. Then the principal eigenvalue
of the operator $-(\M + a(x)) $ in $C_{per}(\R)$ is given by
\begin{align}
\label{variational eigenvalue R} \lambda_1(a) &=
\inf_{\|\varphi\|_{L^2(\R)} = 1} - \int_{\R} (\M[\varphi] + a(x)
\varphi) \varphi
\\
\bet
&= \min_{\varphi \in C_{per}(\R)} - \frac{\int_{-R}^R (\M[\varphi]
+ a(x) \varphi) \varphi}{\int_{-R}^R \varphi^2}.
\end{align}
\end{lemma}
\unskip

\begin{proof}
By Theorem~\ref{cdm.th.pev} we know that there exists
a unique principal eigenvalue $\lambda_1(a) $ of the operator $- (
\M + a)$ in $C_{per}(\R)$. Let $\phi_1 \in C_{per}(\R)$ denote a
positive eigenfunction associated with $\lambda_1(a) $. We
normalize $\phi_1$ such that
\begin{align}
\label{normalization phi1} \int_{-R}^R \phi_1^2 = 2 R.
\end{align}
On the other hand, the quantity
$$
\tilde \lambda_1(a) = \inf_{\varphi \in C_{per}(\R)} -
\frac{\int_{-R}^R (\M[\varphi] + a(x) \varphi)
\varphi}{\int_{-R}^R \varphi^2}
$$
is also an eigenvalue of $-( \M + a ) $ on $C_{per}(\R)$ with a
positive eigenfunction. By uniqueness of the principal eigenvalue,
$\lambda_1(a) = \tilde \lambda_1(a)$.

We claim that
$$
\inf_{\|\varphi\|_{L^2(\R)} = 1} - \int_{\R} (\M[\varphi] + a(x)
\varphi) \varphi \le \lambda_1(a).
$$
Indeed, for $r>0$ let $\eta_r \in C_0^\infty(\R)$ be such that
$0\le\eta_r\le1$, $\eta_r(x) = 1$ for $|x|\le r$, $\eta_r(x)=0$
for $|x|\ge r+1$. It will be sufficient to show that
\begin{align}
\label{convergence quotient}
\lim_{r\to+\infty} \frac{ \int_\R (
\M[\phi_1 \eta_r] + a \phi_1 \eta_r ) \phi_1 \eta_r}{\int_\R
(\phi_1 \eta_r)^2} = -\lambda_1(a).
\end{align}
By \eqref{normalization phi1} we have
\begin{align}
\label{l2 norm} \int_\R (\phi_1 \eta_r)^2 = 2 r + O(1) \quad
\hbox{ as } r\to+\infty.
\end{align}
Let $0<\theta<1$. Then
\begin{align}
\nonumber |\M[\phi_1](x) - \M[\phi_1 \eta_r] | & \le \|\phi_1
\|_{L^\infty} \int_{|x-z|\ge r} |J(z)| \, d z
\\
\bet
\nonumber & \le \|\phi_1 \|_{L^\infty} \int_{|z|\ge(1-\theta)r}
|J(z)| \, d z\quad \hbox{for all } |x|\le \theta r
\\
\bet
\label{integral to 0} &= o(1)\quad \hbox{uniformly for all }
|x|\le \theta r .
\end{align}
We split the integral
\begin{align}
\label{split} \int_\R ( \M[\phi_1 \eta_r] + a \phi_1 \eta_r )
\phi_1 \eta_r = \int_{|x|\le \theta r} \ldots d x + \int_{|x|\ge
\theta r} \ldots d x .
\end{align}
Using $\eta_r(x) = 1 $ for $|x|\le \theta r $ and \eqref{integral
to 0} we see that
\begin{align*}
\int_{|x|\le \theta r} ( \M[\phi_1 \eta_r] + a \phi_1 \eta_r )
\phi_1 \eta_r & = \int_{|x|\le \theta r} ( \M[\phi_1 \eta_r] + a
\phi_1 ) \phi_1
\\
\bet
&= \int_{|x|\le \theta r} ( \M[\phi_1] + a \phi_1 + o(1) ) \phi_1
\\
\bet
&= -2 \theta \lambda_1(a) r + o(r) \quad \hbox{ as } r\to+\infty.
\end{align*}
The second integral in \eqref{split} is bounded by
\begin{align}
\label{rest} \left| \int_{|x|\ge \theta r} ( \M[\phi_1 \eta_r] +
a \phi_1 \eta_r ) \phi_1 \eta_r \right| \le C (1-\theta) r.
\end{align}
Thus from \eqref{l2 norm}--\eqref{rest} we conclude that
\begin{align*}
\left|\frac{ \int_\R ( \M[\phi_1 \eta_r] + a \phi_1 \eta_r )
\phi_1 \eta_r}{\int_\R (\phi_1 \eta_r)^2} + \lambda_1(a) \right|
\le C (1-\theta) + o(1),
\end{align*}
which proves \eqref{convergence quotient}.

\pagebreak

To establish \eqref{variational eigenvalue R} it remains to verify
that
\begin{align}
\label{other ineq eig} \lambda_1(a) \le - \frac{\int_\R
(\M[\varphi] + a(x) \varphi) \varphi}{\int_\R \varphi^2} \quad
\hbox{for all } \varphi \in C_c(\R).
\end{align}
By uniqueness of the principal eigenvalue we have
\begin{align}
\label{def lambda1k} \lambda_1(a) = \inf_{\varphi \in
C_{per}(\Omega_k)} - \frac{\int_{- kR}^{k R} (\M[\varphi] + a(x)
\varphi) \varphi}{\int_{-k R}^{k R} \varphi^2},
\end{align}
where
$$
\Omega_k = (-kR, kR) \quad \hbox{ for $k \ge 1$}
$$
and $C_{per}(\Omega_k)$ is the set of continuous $2 k R$-periodic
functions on $\R$.

Fix $\varphi \in C_c(\R)$ and consider $k$ large enough so that
${\rm supp}(\varphi) \subseteq \Omega_k$. Consider now $\varphi_k$ the
$4 kR$-periodic extension of $\varphi$. Since $\varphi_k \in
C_{per}(\Omega_{2k})$, \eqref{def lambda1k} yields
\begin{align}
\label{lambda1k} \lambda_1(a) \le - \frac{\int_{- 2kR}^{2k R}
(\M[\varphi_k] + a(x) \varphi_k) \varphi_k}{\int_{-2k R}^{2k R}
\varphi_k^2} = - \frac{\int_{\R} (\M[\varphi_k] + a(x) \varphi)
\varphi} {\int_{\R} \varphi^2}.
\end{align}
For $|x|\le kR$ we have
\begin{align*}
| \M[\varphi_k](x) - \M[\varphi](x) | \le \|\varphi\|_{L^\infty}
\int_{|y|\ge 2kR} |J(x-y)| \, d y \le \|\varphi\|_{L^\infty}
\int_{|z|\ge k R} |J(z)|\, d z.
\end{align*}
Hence
\begin{align}
\label{limit periodic extension} \lim_{k\to + \infty} \int_\R
(\M[\varphi_k] + a(x) \varphi) \varphi = \int_\R (\M[\varphi] +
a(x) \varphi) \varphi.
\end{align}
Thanks to \eqref{lambda1k} and \eqref{limit periodic extension}, we
conclude the validity of \eqref{other ineq eig}.
\qquad\end{proof}

\begin{lemma}
\label{convergnce eig}
Assume $J$ satisfies \eqref{J1}, \eqref{J2},
and \eqref{J symmetric} and that $a\colon\R\to\R$ is a
$2R$-periodic, Lipschitz function. Let $\lambda_{r,y}$ be the
principal eigenvalue of \eqref{truncated eigenvalue} for
$$
\Omega_{r,y} = B_r(y)
$$
and let $\lambda_1(a)$ denote the principal eigenvalue of $-(\M +
a(x))$ in $C_{per}(\R)$. Then
$$
\lim_{r\to + \infty} \lambda_{r,y} = \lambda_1(a).
$$
Moreover, the applications $y\mapsto \lambda_{_{r,y}}$ and
$y\mapsto \varphi_{_{r,y}}$ are periodic. The periodicity of the
application $y\mapsto \varphi_{r,y}$ is understood as follows:
$$
\varphi_{_{r,y+2R}}(x)=\varphi_{_{r,y}}(x-2R).
$$
\end{lemma}
\unskip

{\it Proof}.
For convenience we write
$$
\lambda_{r} = \lambda_{r,y}
$$
and let $\varphi_r$ be a positive eigenfunction of
\eqref{truncated eigenvalue} in $\Omega_r$.

By the variational characterization \eqref{variational eigenvalue
Omega} we see that $r \mapsto\lambda_r$ is nonincreasing, and
hence $\lim_{r\to+\infty} \lambda_r$ exists. Moreover, using
\eqref{variational eigenvalue R} we have
\begin{align}
\label{ineq lambdar} \lambda_r \ge \lambda_1(a) \quad \hbox{for
all } r>0.
\end{align}
Let $\phi_1 \in C_{per}(\R)$ be a positive eigenfunction of $-( \M
+ a(x) ) $ with eigenvalue $\lambda_1(a)$ normalized such that
$$
\int_{-R}^R \phi_1^2 = 2 R.
$$
Let $\eta_r \in C_0^\infty(\R)$ be such that $0 \le\eta\le 1$,
\begin{align*}
\hbox{$\eta_r(x) = 1 $ for $|x-y|\le r-1 $,} \qquad
\hbox{$\eta_r(x) = 0 $ for $|x-y|\ge r $}
\end{align*}
and such that $\|\eta_r\|_{C^2(\R)} \le C$ with $C$ independent of
$r$. Arguing in the same way as in the proof of Lemma~\ref{lema
formulas eig R} we obtain
\begin{align*}
\lim_{r\to+\infty} \frac{ \int_\R ( \M[\phi_1 \eta_r] + a \phi_1
\eta_r ) \phi_1 \eta_r}{\int_\R (\phi_1 \eta_r)^2} =
-\lambda_1(a).
\end{align*}
Since
$$
\lambda_r \le - \frac{ \int_\R ( \M[\phi_1 \eta_r] + a \phi_1
\eta_r ) \phi_1 \eta_r}{\int_\R (\phi_1 \eta_r)^2}
$$
we conclude that
$$
\lim_{r\to+\infty} \lambda_r \le \lambda_1(a) .
$$
This and \eqref{ineq lambdar} prove the desired result.

Let us now show the periodicity of the applications $y\mapsto
\lambda_{r,y}$ and $y\mapsto \varphi_{r,y}$. Replace $y$ by $y+2R$
in the above problem \eqref{truncated eigenvalue} and let us
denote by $\lambda_{_{r,y+2R}}$ and $\varphi_{_{r,y+2R}}$
the corresponding principal eigenvalue and the associated positive
eigenfunction:
$$
\opm{\varphi_{_{r,y+2R}}}+a(x)\varphi_{_{r,y+2R}}=-\lambda_{_{r,y+2R}}\varphi_{_{r,y+2R}}
\quad \text{ in }B_r(y+2R).
$$
We take the following normalization:
$$
\int_{\O_{_{r,y+2R}}}\varphi_{_{r,y+2R}}^2(x)\,dx =1.
$$

Let us defined $\psi(x):=\varphi_{_{r,y+2R}}(x+2R)$ for any $x\in B_r(y)$.
A short computation shows that
$$
\opm{\psi}(x) =\opm \varphi_{_{r,y+2R}}(x+2R).
$$
Therefore, using the periodicity of $a(x)$, we have
\begin{align*}
\opm{\psi}(x) +a(x+2R)\psi(x)=\lambda_{_{r,y+2R}}\psi \quad \text{in } B_r(y),\\
\bet \opm{\psi}(x) +a(x)\psi(x)=\lambda_{_{r,y+2R}}\psi \quad
\text{in } B_r(y).
\end{align*}
Thus, $\lambda_{_{r,y+2R}}$ is a principal eigenvalue of the
problem \eqref{truncated eigenvalue} with $\Omega_{r,y}=B_r(y)$.
Hence, by uniqueness of the principal eigenvalue we have
$\lambda_{_{r,y}}=\lambda_{_{r,y+2R}}$ and
$\psi=\gamma\varphi_{_{r,y}}$ for some positive $\gamma$. Using
the normalization, it follows that $\gamma=1$. Therefore,
$\varphi_{_{r,y}}(x)=\varphi_{_{r,y+2R}}(x+2R)$; in other words
$$
\varphi_{_{r,y+2R}}(x)= \varphi_{_{r,y}}(x-2R).
$$
\fdem

\begin{remark}
\label{uniform convergence} \rm
The proof of Lemma~\ref{convergnce
eig} yields the slightly stronger conclusion that the convergence
\begin{align*}
\lim_{r\to+\infty} \lambda_{r,y} = \lambda_1(a)
\end{align*}
is uniform with respect to $y \in \R$, since $\lambda_{r,y}$ is
continuous in $y$.
\end{remark}

{\it Proof of Proposition}~\ref{positive lower bound}.
Let $u\ge 0$
be a bounded solution to \eqref{cdm.eq.noloc} such that
$u\not\equiv 0$. By the strong maximum principle
(Theorem~\ref{strong max p a.e.}) we must have $\inf_K u > 0$ for
compact sets $K \subset \R$.

Given $y\in \R$ and $r>0$ we write $\Omega_{r,y} = (y-r,y+r)$,
$\lambda_{r,y}$ the principal eigenvalue of $-(\M + f_u(x,0))$
with Dirichlet boundary condition in $\Omega_{r,y}$ as in
\eqref{truncated eigenvalue}, and $\varphi_{r,y}$ a positive
Dirichlet eigenfunction normalized so that
$$
\int_{\Omega_{r,y}} \varphi_{r,y}^2 = 1.
$$
Since the principal eigenvalue $\lambda_1 := \lambda_1(f_u(x,0))$
of $-(\M + f_u(x,0))$ with periodic boundary conditions is
negative by hypothesis, by Lemma~\ref{convergnce eig} and
Remark~\ref{uniform convergence} we may fix $r>0$ large enough so
that
$$
\lambda_{r,y} < \lambda_1/2 \quad \text{for all $y \in \R$}.
$$
Note that for $x \in \Omega_{r,y}$,
\begin{align*}
\M[\gamma \varphi_{r,y}] + f(x,\gamma _{r,y} ) & = -\lambda_{r,y}
\gamma \varphi_{r,y} - f_u(x,0) \gamma \varphi_{r,y} + f(x,\gamma
\varphi_{r,y})
\\
\bet
&\ge - \lambda_1 / 2 \gamma \varphi_{r,y} - f_u(x,0) \gamma
\varphi_{r,y} + f(x,\gamma \varphi_{r,y})
\\
\bet
& \ge 0
\end{align*}
if $0\le \gamma \le \gamma_0$ with $\gamma_0$ fixed suitably
small. For $x \not\in \Omega_{y,r}$ we have $\varphi_{y,r}(x) = 0$
and $\M[\varphi_{r,y}] \ge 0$. Thus
\begin{align}
\label{gamma varphi1 subsolution} \M[\gamma \varphi_{r,y}] +
f(x,\gamma \varphi_{r,y} ) \ge 0 \quad \hbox{in $\R$}
\end{align}
for all $0<\gamma<\gamma_0$.

We claim that
\begin{align}
\label{ineq gamma varphi u} \gamma_0 \varphi_{r,y} \le u \quad
\hbox{in $\R$ \quad for all $y \in \R$}.
\end{align}
This proves the proposition because there is a positive constant
$c$ such that $\varphi_{r,y}(y) \ge c$ for all $y \in \R$ since
the application $y\mapsto\varphi_{r,y}$ is periodic and
$\varphi_{r,y}(y)>0$ for any $y\in[-2R,2R]$.

Now, to prove \eqref{ineq gamma varphi u} fix $y\in \R$ and set
$$
\gamma^* = \sup \{ \, \gamma>0 \, / \, \gamma \varphi_{r,y} \le u
\hbox{ in $\R$} \}.
$$
Since $\inf_K u > 0$ for compact sets $K \subset \R$ and
$\varphi_{r,y}$ has compact support we see that $\gamma^*>0$.
Assume that $\gamma^* < \gamma_0$. Then by \eqref{gamma varphi1
subsolution}, $\gamma^* \varphi_{r,y}$ is a subsolution of
\eqref{cdm.eq.noloc} while $u$ is a solution. By the strong
maximum principle (Theorem~\ref{strong max p a.e.}) either
$\gamma^* \varphi_{r,y} \equiv u$ in $\R$ or $\inf _K (u -
\gamma^* \varphi_{r,y}) >0$ for compact sets $K \subset \R$. The
former case is impossible because $u$ is strictly positive, while
the latter case yields a contradiction with the definition of
$\gamma^*$. It follows that $\gamma^* \ge \gamma_0$ as desired.
\fdem 

\section*{Appendix}

In this appendix we give a short proof of Theorem~\ref{strong max
p a.e.}. We assume that $J$ satisfies \eqref{J1}, \eqref{J2}, $c
\in L^\infty(\R)$, and $u \in L^\infty(\R)$ satisfies
\begin{align}
\nonumber u & \le 0 \quad \hbox{a.e. in $\R$},
\\
\bet
\label{ineq eq} \M[u] + c u & \ge 0 \quad \hbox{a.e. in $\R$}.\tag{A.1}
\end{align}
For $\eps>0$ define
$$
u_\eps(x) = \frac{1}{2\eps} \int_{x-\eps}^{x+\eps} u .
$$
Then $u_\eps$ is continuous in $\R$, $u_\eps \le 0 $, and $u_\eps
\to u $ a.e.\@ as $\eps \to 0$. There are two cases:

(1) for any closed interval $I$ one has $\limsup_{\eps \to 0}
\sup_I u_\eps < 0$, or

(2) for some closed interval $I$ one has $\limsup_{\eps \to 0}
\sup_I u_\eps = 0$.

If case (1) occurs, we see that for all closed intervals $I$ we
have ${{\rm ess} \ \sup}_I u<0$. Assume case (2) holds. Let $I$ be a closed
interval and $\eps_n \to 0$ be such that $\lim_{n\to+\infty}
u_{\eps_n}(x_n) = 0$, where $x_n \in I$ is such that $\sup_I
u_{\eps_n} = u_{\eps_n}(x_n)$. Integrating \eqref{ineq eq} from
$x_n - \eps_n $ to $x_n + \eps_n $ and dividing by $2 \eps_n$, we
have
\begin{align*}
J \star u_{\eps_n} (x_n) \ge u_{\eps_n}(x_n) - \frac{1}{2\eps_n}
\int_{x_n - \eps_n}^{x_n + \eps_n} c u .
\end{align*}
But, since $u \le 0$ a.e.,
\begin{align*}
\left| \frac{1}{2\eps_n} \int_{x_n - \eps_n}^{x_n + \eps_n} c u
\right| \le - \|c\|_{L^\infty} u_{\eps_n}(x_n) \to 0.
\end{align*}
Hence
\begin{align*}
\liminf_{n\to+\infty} J \star u_{\eps_n} (x_n) \ge 0.
\end{align*}
We may assume that $x_n \to x \in I$. Then by dominated
convergence,
\begin{align*}
J \star u_{\eps_n} (x_n) = \int_\R J(x_n - y ) u_{\eps_n}(y) \, d
y \to \int_\R J(x - y ) u(y) \, d y.
\end{align*}
This shows that $u = 0$ a.e.\ in $x - {\rm supp}(J)$. Now, for any
$x_1$ in the interior of $x - {\rm supp}(J)$ we have $J \star u (x_1)
\ge 0$, which shows that $u = 0$ a.e.\  in $x - 2 {\rm supp}(J)$, where
$2 {\rm supp}(J) = {\rm supp}(J) + {\rm supp}(J)$. Note that assumption \eqref{J2}
implies that $k \, {\rm supp}(J)$ covers all of $\R$ as $k\to+\infty$,
where $k \, {\rm supp}(J) $ is defined inductively as $(k-1) \, {\rm supp}(J)
+ {\rm supp}(J)$. Repeating the previous argument we deduce that $u=0$
a.e.\@ in $\R$.
\fdem 

\pagebreak

\end{document}